\documentclass[12pt,a4paper]{article}

\usepackage{graphicx}
\usepackage{amsmath}
\usepackage{amssymb}
\usepackage{enumerate}
\usepackage{color}
\usepackage{natbib}
\newcommand{\eref}[1]{(\ref{#1})}
\newcommand{\Var}{\textrm{Var}}

\newcommand{\dd}{\,\textrm{d}}

\newcommand{\NA}{\textrm{NA}}

\newcommand{\Jyvaskyla}{Jyv{\"a}skyl{\"a}}

\title{Optimal design of observational studies: overview and synthesis}
\author{Juha Karvanen$^1$,
Jarno Vanhatalo$^{2}$,
Kari Auranen$^{3}$,\\
Sangita Kulathinal$^4$ and
Samu M{\"a}ntyniemi$^5$\\
~\\
$^1$ \small Department of Mathematics and Statistics, University of Jyvaskyla,\\ \small \Jyvaskyla, Finland, juha.t.karvanen@jyu.fi\\
$^2$ \small Department  of Mathematics and Statistics and\\ \small Department of Biosciences, University of Helsinki,\\ \small Helsinki, Finland\\
$^3$ \small Department of Mathematics and Statistics and \\ \small Department of Clinical Medicine, University of Turku,\\ \small Turku, Finland\\
$^4$ \small National Institute for Health and Welfare, Helsinki, Finland\\
$^5$ \small Department of Environmental Sciences, University of Helsinki,\\ \small Helsinki, Finland 
}
\begin{document}
% \title{Optimal design of observational studies: overview and synthesis}
% %\title{Bayesian design of observational studies: overview and synthesis}
% %\title{Optimal design of observational studies}
% %\title{Optimal design of observational studies: a new field of statistics}
% %\title{Recent advances in optimal design of observational studies}
% \author{Juha Karvanen$^1$,
% Jarno Vanhatalo$^{2,3}$,
% Kari Auranen$^{4}$,\\
% Sangita Kulathinal$^2$ and
% Samu M{\"a}ntyniemi$^6$\\
% ~\\
% $^1$ Department of Mathematics and Statistics, University of Jyvaskyla,\\ \Jyvaskyla, Finland, juha.t.karvanen@jyu.fi\\
% $^2$ Department  of Mathematics and Statistics, University of Helsinki,\\ Helsinki, Finland\\
% $^3$ Department of Biosciences, University of Helsinki, Finland\\
% $^4$ Department of Mathematics and Statistics and \\Department of Clinical Medicine, University of Turku,\\ Turku, Finland\\
% $^5$ National Institute for Health and Welfare, Helsinki, Finland\\
% $^6$ Department of Environmental Sciences, University of Helsinki,\\ Helsinki, Finland }
% %(The first author randomized the order of the other authors.\\ Tell him if you wish to replace this random design by\\ a deterministic design.)}

\maketitle

\begin{abstract}
We review typical design problems encountered in the planning of observational studies and propose a unifying framework that allows us to use the same concepts and notation for different problems.
In the framework, the design is defined as a probability measure in the space of observational processes that determine whether the value of a variable is observed for a specific unit at the given time. The optimal design is then defined, according to Bayesian decision theory, to be the one that maximizes the expected utility related to the design. 
We present examples on the use of the framework and discuss methods for deriving optimal or approximately optimal designs.

\end{abstract}

%We review a wide variety of design problems that can be cast in the proposed framework, including sample size determination, subsample selection, selection for re-measurements, choice of measurement times, and optimization of spatial measurement networks.\\
% ~\\
% AMS 1991 subject classifications. Primary 62C12; secondary 62D05, 62K05\\
% %Key words and phrases: Bayesian design, Decision theory, Sampling, Sequential design, Survey, Value of information\\
% Key words and phrases: Decision theory, Sampling, Sequential design, Survey, Value of information\\

\section{Introduction}
Observational data may arise as the end-product of a well-planned observational study or as the side-product of some operational process. In both cases, data collection can be characterized by three dimensions: the units from which the data are collected, the variables to be measured and the times of the measurement. 
%Unlike in experimental studies, where the analyst has full control over all these dimensions, in observational studies 
The control over these dimensions vary case by case. In a planned observational study, such as a tradional survey, the investigator may decide the units to be included, the variables to be measured and the times of measurement. In operational data collection, the units, variables and times are often determined by the process itself. Transactional data on purchases by customers is an example of this setup. 

The cost of data collection is a major factor restricting the accumulation of scientific knowledge. It is therefore of primary interest that  studies are planned and implemented efficiently. Smart decisions on the study design maximize the expected precision of parameter estimates with the given budget or, alternatively, minimize expected costs while ensuring that the goals of the study are met.

Due to differences in data collection procedures across different application areas, the design of observational studies has evolved without a common framework. This slows down the development and exchange of new ideas and methods.
Hence, a unifying framework for observational studies would enhance the development of new methods by highlighting the similarities and best practices among different applications. 
Here, we propose such a framework and show how it can be applied to wide variety of problems. We also discuss methods to find optimal observational designs and how the proposed framework differs from the widely used optimal experimental design framework \citep{Chaloner:Bayesianexperimental}.

As motivating examples, we consider three problems where the design decisions are related to: (a) the number of units, (b) the individual units to be selected,  and (c) the individual measurement times. We will later show how these problems can be treated under the same framework.

The first example considers collection of fishing data using clustered sampling. Data on the key properties (length, weight, maturity stage and sex) of the fish catch are needed for decisions on fishery policy. The design for the collection of fishing data includes specifying the number of fishing vessels, the number of fishing journeys and the number of individual fish to be measured. There is a cost for each additional vessel, journey and individual to be included in the data collection.

The second example considers subsample selection in medical studies. Specialized medical measurements, such as whole genome sequencing, magnetic resonance imaging or recording maximal oxygen consumption, may be expensive to carry out for a large number of individuals. Therefore an investigator may consider choosing a subsample for the expensive measurements. Choosing a random sample is a natural choice but some other choice often has a better cost-efficiency.

The third example considers the choice of the measurement times in longitudinal studies carried out in medical and social sciences. The characteristics of an individual, such as disease status, health behaviour, socio-economic position and attitudes, change during the course of time.  Measuring the time-varying process very frequently makes sure that all relevant information is gathered but is a poor choice from efficiency perspective if there is a cost associated with each measurement.

We define an observational design as a probability measure in the space of observational processes that determine whether the value of a variable is observed for a specific unit at a given time. This definition covers different types of random sampling as well as designs where the units to be measured are decided in such a way that the expected utility is maximized. While utility functions are commonly applied in experiemental design, they have not been systematically used to guide the design of observational studies.  We formulate the principles of optimal design of observational studies (ODOS) starting from  Bayesian decision theory \citep{Raiffa+Schlaifer,lindley1972bayesian,berger1985statistical}. ODOS can be viewed as a part of a wider decision theoretic framework where the (monetary) value of information \citep{Raiffa+Schlaifer,Lindley:decisions,Eidsvik:VoI} is used as the leading principle to make decisions about collecting additional data. This extends the scope of applications from scientific research to decision making in business and society.

The rest of the paper is organized as follows.
The notation and key concepts are introduced in Section~\ref{sec:notation}. The proposed framework is formulated in Section~\ref{sec:principles}. Computational methods for finding optimal designs are discussed in Section~\ref{sec:finding}.
%Issues related to data analysis are considered in Section~\ref{sec:dataanalysis}.
The motivating examples are presented using the framework in Section~\ref{sec:applications} where other applications are also reviewed. Open problems and future directions are discussed in Section~\ref{sec:discussion}.

\section{Setting and notation} \label{sec:notation}
In this section, we define concepts that can be used to describe a wide variety of design problems, including the three examples given in Introduction. As a starting point, it is required that the high level objectives of the study are known and the study population, study variables and study period are specified. The study population $\mathcal{I}$ is a set of all observational units $i \in \mathcal{I}$ that can be selected for the study. In many applications especially in medicine and social science, the units are individuals and the study population is discrete and finite $\mathcal{I}=\{1,\ldots,N\}$, where $N$ is the size of the population. In other applications, especially in spatial statistics, the study population (study space) is continuous but has a finite area or volume. The units can then be understood to be locations. The study variables, indexed as $j=1,\ldots,J$,  are variables that can be measured (but are not necessarily measured) during the study. The study period is a union of all time intervals during which measurements can be potentially carried out.

Following the idea of the three design dimensions, let $\{x_{ij}(t)\}$, where $i$, $j$ and $t$ refer to unit, variable and time, respectively, denote a stochastic process with values in space $\mathcal{X}$.
%let continuous time process $x_{ij}(t)$ in space $\mathcal{X}$ denote the true underlying value of variable $j$ for unit  $i \in \mathcal{I}$ at time $t$.
The observational process $\{r_{ij}(t)\}$ in space $\mathcal{R}$ is defined as follows
\begin{equation} \label{eq:rdefinition}
r_{ij}(t)=
 \begin{cases}
  1 & \textrm{if } x_{ij} \textrm{ is measured at time } t, \\
  0 & \textrm{otherwise}.
 \end{cases}
\end{equation}
The measurements for $x_{ij}(t)$ can be made only at time points in a finite or infinite set  $\mathcal{T}^{ij}$.  A collection of observational processes $\mathbf{r}=\{\{r_{ij}(t): t \in \mathcal{T}^{ij}$\}: $j \in \{1,\ldots,J\}$, $i \in \mathcal{I}\}$ is called a measurement plan. The measured data can then be defined using the standard notation for missing data
\begin{equation} \label{eq:x*definition}
x_{ij}^*(t)=
 \begin{cases}
  x_{ij}(t) & \textrm{if } r_{ij}(t)=1 \\
  \NA & \textrm{if } r_{ij}(t)=0.
 \end{cases}
\end{equation}
%For a fixed time point $t \in \mathcal{T}^j$, variable $r_{ij}(t)$ is a selection indicator that determines the units for which variable $j$ is measured at time $t$.

An observational design can be informally defined as a strategy for determining the observational processes $\{r_{ij}(t)\}$ for all $i$, $j$ and $t$.
%More formally, design $\eta$ is defined by a probability measure $q_\eta$ on the space of observational processes $\mathcal{R}$.
More formally, a design $\eta$ is defined as a probability measure on the space of observational processes $\mathcal{R}$. This definition differs from the definition of experimental design as a probability measure on space~$\mathcal{X}$ \citep{Chaloner:Bayesianexperimental}. The induced subspace $\mathcal{R}_\eta = \{\mathbf{r} \in \mathcal{R}: \eta(\mathbf{r}) >0\}$ contains all measurement plans that are possible under the design $\eta$.

The above definition allows randomness in the selection of units, variables and measurement times. A design $\eta$ is deterministic if $|\mathcal{R}_\eta|=1$, otherwise it is random. For a deterministic design, the set $\mathcal{R}_\eta$ contains only one measurement plan. For a random design, $\mathcal{R}_\eta$ is a set of measurement plans and  $\eta$ defines the probabilities for each of them to be chosen. As an example,  consider simple random sampling of units for the measurement of variable $j$ at a fixed time point. In this case, design $\eta$ is characterized by sample size $n_\eta$ and selection indicators $r_i$, $i \in \mathcal{I}$, are sufficient to identify a particular measurement plan.
%as the measurement time(s) and the variables to be measured are the same for all units.
The set $\mathcal{R}_\eta$ contains all measurement plans $\mathbf{r}$ for which $\sum_i r_i = n_\eta$. The probability of each measurement plan to be chosen is the same $\eta(\mathbf{r})=1/|\mathcal{R}_\eta|$. A random design is realized when one measurement plan from set $\mathcal{R}_\eta$ is chosen by random sampling. The null design $\eta_0$ is a special case where no new data are collected.

In definition~\eref{eq:rdefinition}, time $t$ refers to the time in the process $x_{ij}(t)$ as well as the time when the measurement is actually made. In some cases, these two times may have to be treated separately. We then define
\begin{equation} \label{eq:rdefinition2}
r_{ij}(t,s)=
 \begin{cases}
  1 & \textrm{if the value of } x_{ij}(t) \textrm{ is measured at time } s, \\
  0 & \textrm{otherwise}.
 \end{cases}
\end{equation}
The retrospective collection of medical history is an example of a situation where the two time axes could be used. Time $s$ refers to  when the individual was interviewed and the medical records were checked while time $t$ pertains to the time of the recorded event.

We allow  data to have a hierarchical structure with levels $k=1, \ldots, K$ so that $N(k)$ is the number of clusters at level $k$ and $k=1$ corresponds to the unit level. The hierarchical structure is specified by an indicator variable $z_{ikl}$ defined as
\begin{equation} \label{eq:clusterindicator}
z_{ikl}=
 \begin{cases}
  1 & \textrm{if unit } i \textrm{ belongs to cluster } l \textrm{ on the level } k, \\
  0 & \textrm{otherwise}.
 \end{cases}
\end{equation}
The hierarchical structure is assumed to remain constant in time but this assumption can be relaxed if needed.
%Is something additional needed for spatial problems?

The general formulation of observational designs allows us to present a wide variety of problems using the same framework. For a hierarchical data structure, the problem may include deciding the sample sizes for each level of the hierarchy as in our first example. There is then only one time point and all variables are measured for all units drawn to the sample. In multi-stage studies, the problem may be the selection of units for the measurement of an expensive variable (second example). The measurement times are considered fixed and a hierarchical structure plays no role. In longitudinal studies, the problem may be the choice of the measurement times (third example). All variables are measured for all units in the study. In replication studies, the problem may be the determination of sample size and selection of variables on the basis of earlier studies. Units are then selected by simple random sampling.

\section{Framework for optimal observational designs} \label{sec:principles}
\subsection{Expected utility of a design}
We assume that prior knowledge on the model parameters $\theta$ is expressed by a distribution $p(\theta)$ defined in parameter space $\Theta$ and that our knowledge on the relationship between the data and the parameters is described by a model $p(\mathbf{x} \vert \theta)$.
The benefits of the study will be evaluated using a utility function and the costs will be measured explicitly in terms of resources spent.

A general framework for planning cost-efficient studies is Bayesian optimal design \citep{Raiffa+Schlaifer,lindley1972bayesian,Chaloner:Bayesianexperimental}, which maximizes the expected utility obtained given the current information about the problem and the cost structure of data collection. This idea has its origins in experimental design but here we adopt it in an observational setting. If data $\mathbf{x}^*_0$ have already been observed, the posterior probability distribution $p(\theta \vert \mathbf{x}^*_0)$ is proportional to  $p(\mathbf{x}^*_0 \vert \theta )p(\theta )$ and describes our current knowledge about the model parameters. In the special case where no previous data $\mathbf{x}^*_0$ have been collected, the current knowledge is described solely by  the prior distribution $p(\theta)$. We consider a design $\eta$ for collecting new data $\mathbf{x}^*_1$ and estimate the expected utility achieved by the design. The utility may depend on the model parameters, the new and the current data, the design and the measurement plan realized under the design for collecting the new data.  If the data are used for decision making, the utility depends also on a decision $d  \in \mathcal{D}$, where $\mathcal{D}$ is the set of possible decisions. Combining these, the utility function can be written in a general form as $U(d,\theta,\mathbf{x}^*_1,\mathbf{r},\eta,\mathbf{x}^*_0)$.
The task is then to find the optimal design, i.e. the one that has the largest maximum expected utility, $\bar{U}(\eta, \mathbf{x}^*_0)$. For any $\eta$, $\bar{U}(\eta, \mathbf{x}^*_0)$ is obtained by conditioning on the existing data $\mathbf{x}^*_0$, marginalizing over measurement plans $\mathbf{r}$ in $\mathcal{R}_\eta$, new data $\mathbf{x}^*_1$ and parameters $\theta$, and maximizing over decision $d$.

In general, the utility function entails both the value of the new information and the costs of data collection. Often, however, it is convenient to define a value function $v(d,\theta,\mathbf{x}^*_1,\mathbf{r},\eta,\mathbf{x}^*_0)$  and a cost function $C(\theta,\mathbf{x}^*_1,\mathbf{r},\eta,\mathbf{x}^*_0)$ separately. In many cases, the value function depends only on the decision and the parameters $\theta$ and can be written as $v(d,\theta)$. The new data $\mathbf{x}^*_1$ increase our knowledge on $\theta$ and thus allow us to make better decisions.
Costs may vary across units and may depend on the existing data $\mathbf{x}^*_0$ or the data $\mathbf{x}^*_1$ yet to  be collected. Often, the costs depend only on the measurement plan, in which case the cost function is simply $C(\mathbf{r})$. Costs are usually defined in terms of money but in some applications it is natural to use time instead of money \citep{karvanen2007experimental,karvanen2009approximate}.

The new data $\mathbf{x}^*_1$ is a random variable in space $\mathcal{X}_1^*$. Their predictive distribution before the data are collected  according to measurement plan $\mathbf{r}$ can be written as
 \begin{equation} \label{eq:newdata}
  p(\mathbf{x}^*_1 \vert \mathbf{r}, \mathbf{x}^*_0) = \int_{\Theta} p(\mathbf{x}^*_1\vert \theta, \mathbf{r},\mathbf{x}^*_0)p(\theta \vert \mathbf{x}^*_0) \dd \theta,
 \end{equation}
 where $p(\mathbf{x}^*_1\vert \theta, \mathbf{r},\mathbf{x}^*_0)$ is the model for the new data. In many applications, $\mathbf{x}^*_1$ and $\mathbf{x}^*_0$ may be assumed to be conditionally independent given $\theta$. However, this does not hold for many time series or spatial problems. On the basis of the new data, our knowledge on the model parameters will be updated to a posterior distribution $p(\theta \vert \mathbf{x}^*_1,\mathbf{x}^*_0, \mathbf{r})=p(\theta \vert \mathbf{x}^*_1,\mathbf{x}^*_0)$, where conditioning on $\mathbf{r}$ is not needed because, according to definition~\eref{eq:x*definition}, the same information can be deducted from $\mathbf{x}^*_1$.

Next we write the maximum expected utility in an explicit form. First we consider deterministic designs and then generalize to random designs.
For deterministic designs, the marginalization over measurement plans is not needed because the design contains only one measurement plan, i.e., $\mathcal{R}_\eta=\{ \mathbf{r}\}$. The optimal deterministic design $\eta^{\textrm{opt}}$ is then equivalent to choosing the measurement plan that maximizes
\begin{align} \label{eq:utility_odos}
& \bar{U}(\eta, \mathbf{x}^*_0) = \bar{U}(\mathbf{r}, \mathbf{x}^*_0) = \nonumber \\
& \int_{\mathcal{X}_1^*} \left[\max_{d \in \mathcal{D}} \int_{\Theta} U(d,\theta,\mathbf{x}^*_1,\mathbf{r},\mathbf{x}^*_0) p(\theta \vert \mathbf{x}^*_1,\mathbf{x}^*_0 ) \dd \theta\right] p(\mathbf{x}^*_1\vert \mathbf{r}, \mathbf{x}^*_0)  \dd \mathbf{x}^*_1
\end{align}
under given constraints. Constraints on decisions and design can be defined either by incorporating them into the utility function or by restricting the set of possible decisions and designs. If $\eta$ is a random design, the expected utility can be expressed as a weighted average over all possible measurement plans
\begin{equation} \label{eq:utility_randomdesign}
 \bar{U}(\eta, \mathbf{x}^*_0) =  \int_{\mathbf{r} \in \mathcal{R}_{\eta}}  \eta(\mathbf{r}) \bar{U}(\mathbf{r} , \mathbf{x}^*_0) ,
 \end{equation}
where $\bar{U}(\mathbf{r} , \mathbf{x}^*_0)$ is defined in equation~\eref{eq:utility_odos}. In a general form, the optimal design 
\begin{equation}\label{eq:optimal_design}
\eta^{\textrm{opt}} = \max_{\eta} \left( \bar{U}(\eta, \mathbf{x}^*_0)  \right),
\end{equation}
is the one that leads to the largest maximum expected utility. In this definition, the utility function entails both the value of the new information and the costs of data collection.

The frequentist approach to ODOS can be presented similarly to criterion~\eref{eq:optimal_design}  with the exception that, instead of the prior and posterior distributions $p(\theta \vert \mathbf{x}^*_0)$  and $p(\theta \vert \mathbf{x}^*_1,\mathbf{x}^*_0)$ in \eqref{eq:utility_odos},  only point estimates of the parameters are considered. If $\hat{\theta}_0$ stands for the initial estimate of $\theta$ based on current data $\mathbf{x}^*_0$, the maximum expected utility of a deterministic design can be written as
\begin{align} \label{eq:utility_odos_freq}
& \bar{U}(\eta, \mathbf{x}^*_0) = \bar{U}(\mathbf{r}, \mathbf{x}^*_0)  = \int_{\mathcal{X}_1^*} \left[ \max_{d \in \mathcal{D}} U(d,\hat{\theta}(\mathbf{x}^*_1,\mathbf{x}^*_0),\mathbf{x}^*_1,\mathbf{r},\mathbf{x}^*_0) \right] p(\mathbf{x}^*_1 \vert \mathbf{r}, \hat{\theta}_0) \dd \mathbf{x}^*_1,
\end{align}
where $\hat{\theta}(\mathbf{x}^*_1,\mathbf{x}^*_0)$ is the point estimate of $\theta$ based on all the data available after the new study. Both the Bayesian~\eref{eq:utility_odos} and frequentist criteria~\eref{eq:utility_odos_freq} require an integration over the future data, and maximization over the decision alternatives. Integration with respect to the prior and the posterior distribution of the model parameters are naturally needed only when applying the Bayesian criterion.

\subsection{Optimal design under cost or utility constraints}

Often the utility function is assumed not to depend on the costs of a design. Typically the reason is that the researchers operate with a fixed budget and the value of the study cannot be easily measured in monetary terms. In this case, one considers separately the expected costs of the design, which in the case of a deterministic design are given by
\begin{align} \label{eq:costs_odos}
\bar{C}(\eta, \mathbf{x}^*_0 ) = \bar{C}(\mathbf{r}, \mathbf{x}^*_0) = & \int_{\mathcal{X}_1^*} \int_{\Theta} C(\theta,\mathbf{x}^*_1,\mathbf{r},\mathbf{x}^*_0) p(\theta \vert \mathbf{x}^*_1,\mathbf{x}^*_0 ) \dd \theta p(\mathbf{x}^*_1\vert \mathbf{r}, \mathbf{x}^*_0)  \dd \mathbf{x}^*_1,
\end{align}
If the costs depend only on the measurement plan, the integration can be omitted and the costs are fixed already before data collection. For the null design, it is natural to assume that $\bar{C}(\eta_0,\mathbf{x}^*_0)=0$.%$C(\eta_0 \vert \mathbf{x}^*_0)=0$.

There are two important special cases of the observational design where costs and benefits are separated:
\begin{enumerate}[(A)]
 \item Maximize the expected utility of the study given a maximum expected cost $C_0$
 \begin{equation*}
  \max_{\eta} \bar{U}(\eta, \mathbf{x}^*_0)  \textrm{ given } \bar{C}(\eta,\mathbf{x}^*_0 ) \leq C_0,
 \end{equation*}
where the utility function is assumed not to depend on the costs. 
\item Minimize the expected cost of the study, given a lower limit $U_0$ for the maximum expected utility
\begin{equation*}
 \min_{\eta} \bar{C}(\eta,\mathbf{x}^*_0 ) \textrm{ given } \bar{U}(\eta, \mathbf{x}^*_0)  \geq U_0,
\end{equation*}
where the utility function is assumed not to depend on the costs.
\end{enumerate}
%Alternatives (A) and (B) are dual problems. Under (A), the researcher has a fixed budget and the research question is determined. A design is selected to give maximal information on the research question, for instance, to minimize the variance of the parameters of interest.

Both of these problems could be formulated as maximization of a general utility function \eqref{eq:optimal_design}. For example, problem (A) is the same as maximizing a utility function that is a step function with respect to design costs so that the utility drops to minus infinity if the costs exceed the limit $C_0$ entailing an impossible design. 
However, often these problems are more intuitive and easier to deal with when defined as above.

\subsection{Value of information analysis}
%{\color{red} We could actually define the value of information analysis as one "important" special case of the general design of experiment problem. Hence, this could be defined as problem (C). I have now modified the text accordingly}

The general form of utility function, which includes also the costs of design, is typical in business and societal decision making, where both the outcome following the decision and the costs of the design can be measured in monetary terms. In this case, the design problem can alternatively be formulated as maximization of the (monetary) value of information \citep{Raiffa+Schlaifer,Lindley:decisions,Eidsvik:VoI}. We illustrate this with a deterministic measurement plan. First, a value function $v(d,\theta,\mathbf{x}^*_1,\mathbf{x}^*_0)$ is defined as the monetary value of decision $d$ for each combination of the parameter values and data obtained through measurement plan $\mathbf{r}$.
The costs of data collection, that is $C(\theta,\mathbf{x}^*_1,\mathbf{r},\mathbf{x}^*_0)$, are assumed to be additive to the value of the decision. The utility can then be defined to be a function of the difference between the value of a decision and the costs of the design, that is $U(d,\theta,\mathbf{x}^*_1,\mathbf{r},\mathbf{x}^*_0) = U\left(v(d,\theta,\mathbf{x}^*_1,\mathbf{x}^*_0)-C(\theta,\mathbf{x}^*_1,\mathbf{r},\mathbf{x}^*_0) \right)$. Further, the value of information of a measurement plan $\mathbf{r}$ is the maximum expected monetary value of observing its outcome, that is, the price $V$ such that \citep{Lindley:decisions,Eidsvik:VoI}
\begin{align}
\int_{\mathcal{X}_1^*} \left[ \max_{d \in \mathcal{D}} \int_{\Theta} U\left(v(d,\theta,\mathbf{x}^*_1,\mathbf{x}^*_0) - V\right) p(\theta \vert \mathbf{x}^*_1,\mathbf{x}^*_0 )\dd \theta \right]p(\mathbf{x}^*_1\vert \mathbf{r}, \mathbf{x}^*_0)  \dd \mathbf{x}^*_1\nonumber \\
= \max_{d \in \mathcal{D}} \int_{\Theta} U \left( v(d,\theta,\mathbf{x}^*_0) \right) p(\theta \vert \mathbf{x}^*_0 ) \dd \theta.\label{eq:VOI}
\end{align}
Hence, we can define another important special case of observational design:
\begin{enumerate}[(C)]
 \item Maximize the value of information of the study 
 \begin{equation*}
  \max_{\eta} V(\eta, \mathbf{x}^*_0)  
 \end{equation*}
where the value of information $V(\eta, \mathbf{x}^*_0)$ is defined implicitly through \eqref{eq:VOI}. 
\end{enumerate}

If the utility is a linear function of value, $U(v) = a+bv$, the value of information for a measurement plan $\mathbf{r}$ reduces to \citep{Raiffa+Schlaifer}
\begin{equation} \label{eq:voi}
 \bar{v}(\mathbf{r}, \mathbf{x}^*_0)  - \bar{v}(\eta_0, \mathbf{x}^*_0) ,
\end{equation}
where
\begin{equation*}
 \bar{v}(\mathbf{r}, \mathbf{x}^*_0) = \int_{\mathcal{X}_1^*} \left[ \max_{d \in \mathcal{D}} \int_{\Theta} v(d,\theta,\mathbf{x}^*_1,\mathbf{x}^*_0)  p(\theta \vert \mathbf{x}^*_1,\mathbf{x}^*_0 )\dd \theta \right]p(\mathbf{x}^*_1\vert \mathbf{r}, \mathbf{x}^*_0)  \dd \mathbf{x}^*_1.
\end{equation*}
The latter form of the value of information quantifies the increase in the expected value due to using design $\eta$ when the reference level is set by the null design. The value of information is the maximum price that a rational decision maker should be willing to pay for having access to new data collected by experimental plan $\mathbf{r}$ before choosing decision $d$. Hence, the decision maker should compare the value of information to the expected cost of employing the design to decide whether to conduct the data collection or not. For an eligible design, the value of information should be greater than the expected cost of data collection. The null design is the optimal design if no other design is eligible. To extend the value of information analysis for random designs, $\eta$, we would need to take the expectation also over the possible measurement plans.

\subsection{On utility functions}

In principle the utility function should be defined for each study based on the study objectives and constraints. 
However, defining utility function for a specific application is not straightforward in general and most of the optimal design literature concentrates on general purpose utility functions. For example, utility functions developed in Bayesian experimental design can be applied also in ODOS. As an example, we mention few common ones here.

 \citet{lindley1956onthemeasure} proposed using the expected Shannon information \citep{Shannon1948} of the posterior distribution of the model parameters as the utility function. In normal linear regression models, this choice leads to Bayesian D-optimal designs where the utility function depends on the determinant of the expected Fisher information. The same happens in non-linear models if the posterior is approximated by the normal distribution \citep{Chaloner:Bayesianexperimental}. In ODOS, Fisher information can be written as $I_{\mathbf{x}_1^*}(\theta,\mathbf{x}^*_0)$ where the notation $I_{\mathbf{x}^*_1}$ indicates that the expectation is taken with respect to the distribution of the new data. This means that the utility is actually a combination of the expected information of $\mathbf{x}_1^*$ and the observed information of $\mathbf{x}_0^*$. Specifically, the utility function for observational Bayesian D-optimality can be written as
\begin{equation} \label{eq:bayesDoptimality}
 U(\theta,\mathbf{x}^*_1,\mathbf{x}^*_0) = \log \det\left(I_{\mathbf{x}^*_1}(\theta)+I^{\textrm{obs}}(\theta,\mathbf{x}^*_0)\right)
\end{equation}
where $I^{\textrm{obs}}$ stands for the observed information. A-optimality and other alphabetical design criteria can be applied in a similar manner.

\subsection{Missing data}

Optimized designs often lead to data missing by design \citep{wacholder1996case,le1997bayesian}. As the decision whether to measure a variable or not may depend on the observed data (and the prior), the data are missing at random (MAR). This situation is thus different from notorious self-selection \citep{keiding2016perils} where the inclusion probability may depend on unknown variables. As the data are MAR, the missing data mechanism, that is, the design, may be ignored in direct-likelihood and Bayesian inference \citep{Rubin:inferenceandmissingdata}. In particular, methods for handling missing data under the MAR assumption can be applied \citep{littlerubin2002}. The inference is model-based, i.e., averages and other descriptive statistics must not be calculated directly from the sample but should be estimated using the model.

%For clarity of presentation, the framework was defined without unintended missing data. 
In addition to data missing by design, some observations may be missing unintentionally. It would be possible to generalize the framework to include also unintended missing data. Adopting the Bayesian paradigm, all unobservables would be  considered jointly with the model parameters $\theta$. Moreover, under ODOS data missing unintentionally can be taken into account also when planning the design. We can define the model for new data $p(\mathbf{x}^*_1\vert \theta, \mathbf{r},\mathbf{x}^*_0)$ to account for the possibility of not observing $x_{ij}(t)$ even if $r_{ij}(t)=1$. 
%The integration over the unobservables in the utility functions would complicate numerical calculations.

\section{Finding optimal designs}  \label{sec:finding}
In most cases, optimal observational designs cannot be found analytically because the integrals in criterion~\eref{eq:utility_odos} cannot be expressed in a closed form. The optimal design then needs to be searched for using numerical methods. However, numerical optimization over the design space in criterion~\eref{eq:utility_odos}  is often computationally intractable. Thus, ``optimal designs'' are in practice only approximately optimal or the best designs obtained under the search strategy used. Despite this shortcoming, the obtained designs may still be highly cost-efficient as compared with standard (non-optimized) designs.

 \citet{ryan2015review} provides a good overview of the available computational algorithms for Bayesian optimal design. The algorithms are presented in the optimal experimental design context but the general ideas are usually applicable also to observational studies.  Computational methods are needed for the estimation of posterior distributions, for the estimation or approximation of utility functions and for searching over the design space. The first two computational problems are shared between all Bayesian optimal design problems.
The search strategies for ODOS differ from those of experimental design as the design space is often formed by a limited number of units. For deterministic designs, we name two basic search strategies as direct search and design search.

The direct search \citep{karvanen2009optimal,reinikainen2014optimal} is related to the search for exact (discrete) optimal experimental designs and uses heuristic methods, such as the greedy method \citep{dykstra1971augmentation} and modified Fedorov method \citep{cook1980comparison}. In the greedy method applied to problem (A), units are selected one by one as long as the budget allows. The original maximization problem is converted to sequential maximization of designs where only one unit is selected. The approach resembles sequential design with the difference that the measurements are envisaged to be made only after all units have been selected, not directly after each selection as in sequential design. In the iterative replacement method applied to problem (A), the search starts with an initial design which can be obtained e.g. by the greedy method. The selected units  are  then considered  one by  one  and  replaced  by  another  unit  if  that  increases the  value  of the utility function. The procedure is iterated until convergence.
%Could the coordinate-exchange algorithm \citep{meyer1995coordinate} be used?

%Design search is applicable in problems where units are selected to be measured and current data $\mathbf{x}_0^*$ are available.
The design search \citep{drovandi2015principled} aims to first find the optimal design for the corresponding experimental setup. The problem of optimal observational design is replaced by a problem of optimal experimental design where the data $\mathbf{x}_0$ represent experimental variables that can be controlled by the researcher.  The solution for this problem is an optimal experimental design characterized by measurement points $\mathbf{x}_0^{\textrm{opt}}$ and their weights. After this, the actual ODOS problem is considered. The observational design is obtained by selecting the units that minimize the distance between the values in $\mathbf{x}_0^*$ and the optimal experimental design $\mathbf{x}_0^{\textrm{opt}}$. For instance, if the optimal experimental design is a two-point design, half of the observational units are selected from the vicinity of each design point. As each unit can be selected only once, the observational design will differ from the optimal experimental design. The approach is expected to work well if the number of units to be selected is small, the number of units available is large, and the number of variables is small to avoid the curse of dimensionality in the distance measure. Otherwise it may happen that all candidate units are far from the optimal experimental design.

\citet{muller1999simulation} and \citet{muller2004optimal} proposed to use Markov chain sampling to locate the optimal design. They introduced an artificial probability density by normalizing the expected utility as a function of the design and the model parameters and constructed a Markov chain to sample from this distribution. After the convergence of the chain, the optimal design can be identified. In spatial statistics, the design is typically not formally optimized over the full design space but the best design is chosen from a set of intuitively sensible designs such as various space filling designs  \citep{muller2007collecting,Chipeta2017,Chipeta2016} or designs that result from reducing the size of an existing monitoring network \cite{Diggle2006}.

Computational methods are needed also for the estimation of posterior distributions. Well-known computational methods, Markov chain Monte Carlo (MCMC) \citep{muller1999simulation}, sequential Monte Carlo \citep{drovandi2014sequential}, importance sampling \citep{cook2008optimal}, Laplace approximation \citep{ryan2015fully} and approximate Bayesian computation (ABC) \citep{drovandi2013bayesian}, have been used in the estimation of posterior distributions and the expected utility.

The computational requirements are usually lighter for random designs than for deterministic designs. Although the analytical expressions become more complicated  under random designs,
it is possible to replace the enumeration of all possible measurement plans in set $\mathcal{R}_{\eta}$  by random sampling from the set.

\section{Examples and applications} \label{sec:applications}
 In subsections \ref{subsec:samplesize}, \ref{subsec:subsample}  and \ref{subsec:times} we apply the proposed framework to three motivating examples given in Introduction. We then present additional applications related to epidemiology, fishery management and spatial measurement networks.

\subsection{Sample size determination} \label{subsec:samplesize}
Sample size determination \citep{Lindley1997} is a special case where the design choice concerns the number of units to be measured. The problem has been studied by many authors in settings where available resources need to be allocated between either multiple stages or multiple subgroups. Closed-form solutions are often available in the frequentist approach. Model-based survey sampling \citep{chambers2012introduction} can be seen as a design problem where the choices concern sample sizes.

As a simple example, assume that the objective is to estimate a population mean $\theta$ with a given precision from a simple random sample. The utility function can depend, for instance, on the precision of the estimated population mean. There is only one time point and all variables are measured for the units in the sample. If no previous data $\mathbf{x}^*_0$ are available and  the precision is measured by the variance, the utility function is $U(d,\theta,\mathbf{x}^*_1,\mathbf{r},\mathbf{x}^*_0)=-\Var(\theta)$. The design is random and characterized by sample size $N_\eta$. The expected utility $\bar{U}(\eta)$ defined by equations \eref{eq:utility_randomdesign} and \eref{eq:utility_odos} can be then calculated by forward simulation in the following steps:
 \begin{enumerate}
  \item Draw a realization $\theta'$ from the prior $p(\theta)$. 
  \item Generate data $\mathbf{x}^*_1$ of $N_\eta$ observations from the model $p(\mathbf{x}^*_1\vert \theta')$.  
\item Estimate the posterior $p(\theta \vert \mathbf{x}^*_1)$ and record the posterior variance.  
\item Repeat steps 1--3, calculate the average posterior variance $\Var(\theta)$ and use it as an estimate of $\bar{U}(\eta)$.
 \end{enumerate}
The smallest sample size with the required  average posterior variance can be found by applying the binary search or other iterative algorithms.

In the case of clustered data, the problem is to decide the sample sizes  $N(k)$ for hierarchy levels $k=1,\ldots,K$ when the total cost of sampling is given. The design is characterized by sample sizes $N(1),\ldots,N(K)$ and the individual measurement plans by their selection indicators $r_i$, $i=1,\ldots,N(1)$. If the cost per cluster at the $k$th cluster level is $c_k$, the total cost of the design will be
\begin{equation*}
 \sum_{k=1}^K c_k N(k)=\sum_{k=1}^K \sum_{l=1}^{L_k} c_k \mathbb{I}\left( \left(\sum_{i=1}^{N(1)} r_i z_{ikl} \right) > 0 \right),
\end{equation*}
where the indicator variable $z_{ikl}$ is defined in equation~\eref{eq:clusterindicator}, $\mathbb{I}$ denotes an indicator function and $L_k$ denotes the number of clusters at hierarchy level $k$. The expected utility can be estimated by means of simulation as above. Finding the optimal sample sizes for each hierarchy level may pose a computational challenge because of large number of possible combinations.

\citet{tokola2014design} presented an application to clustered fishing data which was briefly described in Introduction. The aim was to collect data on the key properties (length, weight, maturity stage and sex) of the fish catch needed for decisions on the fishery policy. The hierarchy had three levels: fishing vessels, fishing journeys and individual fish. Realistic marginal costs (in euro) were set for each additional vessel, journey and individual to be included in the data collection. There were natural upper limits for the number of vessels and the number of journeys. The objective was to minimize the total cost of data collection when the precision targets for the key properties were derived on the basis of the regulations by European Union (a type~B problem). The authors utilized geometric programming \citep{boyd2007tutorial} and ended up with an interesting solution where all available vessels  should be used but only nine fish per journey should be measured.

Many examples on sample size determination in the frequentist framework can be found in the literature.
\citet{reilly1996optimal} considered optimal sampling strategies for epidemiological two-stage studies.  \citet{mcnamee2002optimal} and \citet{wruck2006sequential} studied optimal designs of two-stage studies for the estimation of the sensitivity and specificity of a diagnostic test. \citet{bekmetjev2012cost} proposed a cost-efficient resampling design for the situation where a fraction of the sample is classified by the same imperfect method twice. \citet{rezagholi2010cost} reviewed cost-efficient designs from the viewpoint of occupational exposure assessment. They considered multi-stage studies and comparisons of measurement methods, concluding that the reviewed studies had used simplified analytical tools. \citet{Sutton:evidencebased} studied sample size selection in a meta-analytic framework and \citet{nikolakopoulou2015planning} extended the approach to network meta-analysis. \citet{tokola2011power} considered an application to the design of a health coaching study. The problem was to select the number of coaches and the number of subjects per coach in such a way that the power of the study is maximized given the total cost.

The optimal sample size determination is studied extensively also in the Bayesian context. For example, \citet{Adcock2017} provides an early analysis of required sample size to obtain desired accuracy for the multinomial ratios, \citet{Joseph1995} examine sample sizes that give highest posterior probability intervals for binomial proportions with a fixed minimum probability and \citet{Rahme2000} studied Bayesian sample size determination when the disease prevalence in a population needs to be estimated with a non-perfect diagnostic test. These examples would fall into category (A) in ODOS if the authors had used proper utility functions but, as shown by \citet{Lindley1997}, the treatment of \citet{Joseph1995} (and also \citet{Rahme2000} in that respect) is not coherent from the decision theoretic point of view and hence could not be formulated under ODOS as such. However, \citet{Lindley1997} extends their example to coherent decision theoretic approach where the coverage of the highest posterior interval, its length and the cost of sampling are tied together with a utility function. After this an optimal decision concerning these three parameters is searched by maximizing the expected utility and the sample size fitting into category (C) in ODOS. The treatment by \citet{Lindley1997} illustrates how design of observational studies is sometimes conducted ``informally'', without defining utility function and, hence, not fulfilling the formal rules or coherence of the decision theoretic ODOS framework.  Nevertheless, with slight reformulation of the problem proper treatment is possible.

\subsection{Subsample selection} \label{subsec:subsample}
%The example presented in Section~\ref{sec:motivatingsubsample} considered optimal subsample in a two-stage study. 
In the basic two-stage setting briefly decribed in our second example in Introduction, variables $\mathbf{x}^*_0 = \{ x_{ij} \}$, $j=1,\ldots,J_0$, $i=1,\ldots,n$ have been measured for the complete sample in the first stage. At the second stage, the problem is to select optimally a subsample of size $n_1<n$ for which variables $x_{ij}$, $j=J_0+1,\ldots,J$ will be measured. Often, the variables at the second stage are expensive to measure and therefore only a subsample can be considered. For instance, a surrogate endpoint may have been measured at the first stage and a small subsample is selected at the second stage to validate the surrogate measurements.

A measurement plan can now be characterized by selection indicators $\mathbf{r}=(r_1,\ldots,r_n)$, where $r_i=1$ if the variables $x_{ij}$, $j=J_0+1,\ldots,J$ are to be measured for unit $i$. The data collected at the second stage is denoted by $\mathbf{x}^*_1 = \{ x_{ij}^* \}$, $j=J_0+1,\ldots,J$, $i=1,\ldots,n$, where $x_{ij}^*=\NA$ if $r_i=0$. If a deterministic design is used, the design problem of type~A takes the form
\begin{align} \label{eq:Asubsampleselection}
 &\max_{\mathbf{r}} \int_{\mathcal{X}_1^*} \left[ \int_{\Theta} U(\theta,\mathbf{x}^*_1,\mathbf{x}^*_0) p(\theta\vert \mathbf{x}^*_1,\mathbf{x}^*_0) \dd \theta  \right] p(\mathbf{x}^*_1 \vert \mathbf{r},\mathbf{x}^*_0)   \dd \mathbf{x}^*_1 \nonumber \\
 & \textrm{given the constraint } \sum_{i=1}^n r_{i} =n_1.
\end{align}

\citet{karvanen2009optimal} studied the selection of individuals for genotyping in a case-cohort study. The variables available at the time of the selection included phenotypic covariates measured at baseline and a survival outcome measured at the end of the follow-up. Blood samples taken at the baseline had been stored so that genotyping did not require the individual to be alive. A frequentist approach with the D-optimality criterion was used. In simulation comparisons, the designs were ranked as expected: D-optimal design, case-cohort design and simple random sampling. However, extreme selection, i.e., selecting individuals with extreme values of phenotypic covariates, worked almost as well as the (approximate) D-optimal design and was recommended as a practical choice. Of note, extreme selection has been a popular design in genetics \citep{lander1989mapping,allison1998extreme,macgregor2006use}. %\citep{lander1989mapping,carey1991linkage,darvasi1992selective,allison1998extreme,van2000power,mcelroy2006comparison,macgregor2006use}.

In addition to medical and epidemiological studies, subsample selection has applications in business. Companies may want to select a subset from their customer database (``big data'') for an experiment or a survey. The selection can be formulated as an optimal design problem. \citet{drovandi2015principled} proposed the optimal selection of a subset to be analyzed as an alternative for parallel computing in the big data context. As examples, they considered datasets on cancer patients, (simulated) mortgage defaults and accelerometer measurements.

\subsection{Choice of measurement times} \label{subsec:times}
Markov models are often used to describe the dynamics of phenomena. Let $X(t)\in \mathcal{X}$  be  a Markov process in a discrete state space $\mathcal{X}= \lbrace{1,\ldots,G}\rbrace$, governed by transition intensities $\theta$ and indexed by continuous time $t$. When the current state of the process is observed only at discrete time points for each study subject, the question emerges about the optimal timing of the measurements as briefly decribed in our third example in Introduction. For illustration, consider measuring the status of bacterial colonisation \citep{mehtala2015optimalSMMR,mehtala2015optimalJRSSC} or the presence of a parasitic infection \citep{nagelkerke1990estimation}. 

In practice, measurements can only be made at discrete time points $\mathcal{T} = \lbrace t_0=0,\ldots,t_M\rbrace$ within a (possibly very large) sampling frame including $n$ individuals. The measurement plan is now characterised by selection indicators ${\mathbf r} = ({\mathbf r}_1,\ldots,{\mathbf r}_n)$, where ${\mathbf r}_{i} = (r_{i}(1),\ldots,r_{i}(M))$ and $r_{i}(m)= 1$ if individual $i$ is included in the sample and the state of his/her process is measured at time $t_m\in \mathcal{T}$, and 0 otherwise. Note that this formulation views also set $\mathcal{T}$ as part of the sampling frame since measurements may only be made in a subset of $\mathcal{T}$.

In the frequentist framework, the optimal design is found as the solution to following problem
\begin{equation*}
\max_{\mathbf{r}} \int U(\hat\theta({\mathbf x}_1^*),{\mathbf x}_1^*)p({\mathbf x}_1^*\vert {\mathbf r},\hat\theta_0)d {\mathbf x}_1^*,
\end{equation*}
with the constraint $\sum_{i=1}^n\sum_{m=1}^M r_{i}(m) = C_0$. Here $C_0$ is the total number of observations that can be made due to practical and/or budget limitations. Often the same number of observations is to be taken from all individuals included in the sample, leading to an additional constraint that for any individual $i$ {\sl included} in the sample $n_1$ measurements are made, i.e., $\sum_{m=1}^M r_{i}(m) = n_1$. The optimal solution will then determine both $n_1$ and the number of individuals sampled as $N = C_0/n_1$.

\citet{mehtala2015optimalSMMR} investigated the optimal sampling interval in a two-state Markov process with equidistant time spacings between consecutive discrete-time observations. The optimal choice of $N$ vs. $n_1$ (cf. above)  was shown to depend on the distribution of the initial condition. In addition, applying a
two-stage design with a utility function of type equation~\eref{eq:bayesDoptimality}, the optimal split of a follow-up study into two stages was investigated. The problem of optimal sampling
times was addressed for a model with more than two states by \citet{mehtala2015optimalJRSSC} in a Bayesian framework.

 \citet{ji2017optimal} considered the optimal timing of measurements in longitudinal data collection when the prior information comes from a pilot study where either densely measured functional data or few randomly scheduled longitudinal measurements are available. Their application was related to the estimation of body-mass index trajectories and the prediction of systolic blood pressure at old age. \citet{cook2008optimal} optimized the measurement times for epidemic processes. \citet{Varis} optimized the monitoring strategy (monthly, biweekly or weekly monitoring) in a lake management problem. \citet{ryan2014towards} and \citet{ryan2015fully} and applied Bayesian optimal design in the planning of the measurement times in memory retention tests and pharmacokinetic studies. Although \citeauthor{ryan2014towards} talk about experimental design, the problem of deciding the measurement times has observational nature and can thus be viewed under the umbrella of ODOS.

\subsection{Selection for re-measurements}
In the basic longitudinal setting, measurements are made at pre-specified time points $\mathcal{T}$ for all units and all variables. We consider a variant where only an optimally selected subset of units will be measured. Let $n_1<n$ be the size of the subset to be measured at time point $t_1$. The data available before time $t_1$ can be expressed as $\mathbf{x}^*_0 = \{ x_{ij}^*(t): t \in \mathcal{T}, t< t_1\}$,  $j=1,\ldots,J_0$, $i=1,\ldots,n$. Our notation defines each re-measurement as a new variable and assumes that the indices of the variables are ordered by the measurement time. A measurement plan can be characterized by selection indicators $\mathbf{r}(t_1)=(r_1(t_1),\ldots,r_n(t_1))$, where $r_i(t_1)=1$ if the variables $j=J_0+1,\ldots,J_1$ are to be measured for unit $i$ at time $t_1$. Writing $\mathbf{r}=\mathbf{r}(t_1)$, the optimization problem becomes similar to subsample selection problem~\eref{eq:Asubsampleselection}. Before the next re-measurement the same problem is considered again but with updated data.

\citet{reinikainen2014optimal} and \citet{reinikainen2016bayesian} considered optimal selection of individuals for longitudinal covariate measurements in epidemiological follow-up studies. The interest was to estimate risk factors associated with time-to-event outcomes. Both frequentist \citep{reinikainen2014optimal} and approximate Bayesian \citep{reinikainen2016bayesian} approaches were applied using the D-optimality criterion. Deaths and drop-outs posed special restrictions on the design. The selection preferred older individuals (higher risk of an event) and individuals with extreme covariate values. A re-analysis of existing epidemiological data
%In the re-analysis of an existing epidemiological data, it was
demonstrated that the sampling size with ODOS can be 10--25\% smaller than under simple random sampling.

\subsection{Spatial problems}
In spatial problems, the measurements are spatially indexed so that $x_{ij}(t)=x_{j}(\boldsymbol{\xi}_i,t)$, where $\boldsymbol{\xi}_i$ is the spatial co-ordinate of the $i$'th spatial observation unit, and the interest typically lies  in spatially indexed latent variables denoted here by $\phi(\boldsymbol{\xi})$. For example, in spatial epidemiology the observations are disease incidences while the latent variables could describe the relative risk for the disease in question \citep{elliot2001,Vanhatalo2010a,Chipeta2016}. The spatial domain is typically either divided into a discrete set of spatial areas, in which case there is a finite number of latent variables and possible observational units, or treated continuously, in which case the latent variables are realizations of a continuous latent process and there is an infinite (uncountable) number of possible observational units within a limited study region. Spatial epidemiology, where the spatial domain would consist of administrative regions, is a typical example of the former and environmental applications of the latter. In both cases we can extend the model parameters to include the latent variables, $\phi\subset \theta$.

Designing observational networks in spatial domain is a long standing subject in spatial statistics. Traditional approach is to construct a space filling design which allows efficient interpolation of the latent field over the study domain \citep{Diggle2006,muller2007collecting,Chipeta2016}. These can be constructed by randomly spreading the sampling locations over the domain or with algorithms producing quasi-random numbers, such as Sobol sequences, that fill the space more uniformly. 
In some applications the interest may be also in the parameters of the spatial field (e.g., covariance function parameters of a spatial Gaussian process) in which case some observation locations should be clustered near each others \citep{muller2004optimal,Chipeta2017}. 
Common to traditional spatial design approaches is that they have proceeded by construction of alternative deterministic or random spatial point patterns which are then compared with a design criteria, the utility function in ODOS framework. A typical utility function is the averaged prediction variance over a spatial area $\mathcal{A}$
\begin{equation*}
U(\theta,\mathbf{x}^*_1,\mathbf{r},\eta,\mathbf{x}^*_0) = \int_{\boldsymbol{\xi}\in \mathcal{A}} \text{Var}(\phi(\boldsymbol{\xi})\vert\mathbf{x}^*_1,\mathbf{x}^*_0)\dd \boldsymbol{\xi}
\end{equation*}
where $\text{Var}(\phi(\boldsymbol{\xi})\vert\mathbf{x}^*_1,\mathbf{x}^*_0)$ is the posterior variance of the latent process. For example, \citep{Diggle2006} compare alternative lattice designs with additional infill points for setting the locations for a fixed number of water monitoring stations in the Baltic sea (a type~A problem). \citet{Chipeta2016} and \citet{Chipeta2017} compare  alternative inhibiroty spatial designs for monitoring the spread of malaria. \citet{Chipeta2016} consider additionally adaptively tuning these designs as data is collected. 

Another typical example is a spatial design where existing monitoring network is decreased by choosing only a subset of stations to further operation \citep{Sanso1999,Diggle2006}. Here the observation stations are fixed number of units and the design decision concerns which of them to retain. \citet{Sanso1999} used a utility function that depends on the accuracy of the rainfall predictions and the cost induced by the number of stations in the network (a type~C problem). In their application the spatial domain was discretized since only a fixed number of stations in predetermined possible locations was considered. Other examples on optimal design of observation networks are provided by, e.g., \citet{muller2004optimal} and \citet{muller2007collecting}.

\subsection{Choice of variables to be studied}
The choice of variables as a design problem occurs at least in two scenarios. In the first one, the objective is to obtain precise predictions at the unit level when there are several potential predictors that could be measured in a new study.  The cost of measurement differs across predictors.
%The problem is to select an optimal set of predictors to be measured.
The second scenario is related to accumulation of scientific information. Researchers should focus their efforts so that the expected  scientific impact of a new study is maximized. This has direct implications not only to decisions on the sample size but also on the variables to be measured in the new study. In both scenarios, a measurement plan can be characterized by selection indicators $\mathbf{r}=(r_1,\ldots,r_J)$, where $r_j=1$ if the variable $j$ is measured for units in the sample. It is possible to extend the problem so that variables to be measured are decided separately for each unit. In this case, a measurement plan is characterized by a selection matrix $\{r_{ij} \}$.

 %It is possible to simultaneously optimize the sample size of the new study.

For the second scenario, \citet{impact} proposed a formal approach to decide which covariates should be measured in a new study in a meta-analytic framework. They compared decision criteria based on conditional power, change of the p-value, change in the lower confidence limit, Kullback-Leibler divergence, Bayes factors, Bayesian false discovery rate or difference between prior and posterior expectations. As an illustration, they considered covariate prioritization based on the results of an existing meta-analysis of genome-wide association studies and made suggestions on the genes to be studied further.

%Karvanen and Sillanp{\"a}{\"a}, a design problem of type C
\subsection{Evaluating value of information}
The value of information can be used to find the maximum price the investigator should pay for additional data. \citet{Eidsvik:VoI} provide a detailed introduction to the theory of value of information in the context of spatial decision analysis. They also address several applications related to e.g. mining, forestry and oil drilling problems. Other examples of value of information analysis can be found from environmental management and health care management \citep{Mantyniemi:Voi, McDonald+Smith,Yokota+Thompson:A,Yokota+Thompson:B}. %and also in the context of earth sciences \citep{Eidsvik:VoI}. 
\citet{Mantyniemi:Voi} analysed the value of information of resolving uncertainty about a biological hypothesis related to the population dynamics of North sea herring. Two alternative hypotheses were considered: the amount of offspring either approaches an asymptotic carrying capacity as the number of parents increases (Beverton-Holt model \citep{BevertonHolt1957}) or peaks at a certain parental population size and then declines (Ricker model \citep{ricker1954stock}). As a result, the maximum price that the fishing industry should be willing to pay for a research that would completely remove the uncertainty about this population dynamic hypothesis was found to be 240 million Norwegian crowns. However, this price was only a 1.6\% increase in the expected utility, compared with optimal solution under the prevailing uncertainty. The authors also invoke the concept of ``price of overconfidence'', the expected loss that would result if ignoring the prevailing uncertainty and optimizing the fishing pressure by assuming that the correct model structure is known.

\section{Discussion} \label{sec:discussion}
We have presented a unifying framework for the design of observational studies under the Bayesian paradigm. In this framework, the observational design is defined as a probability measure in the space of observational processes (measurement plans) which determine whether the values of specific variables are to be observed for specific units at given times. The optimal observational design is defined as the sampling strategy that maximizes the expectation of a pre-defined utility function. We reviewed methods that can be used to find the optimal design and coined the terms direct search and design search. The framework can be used to describe a wide variety of design problems as demonstrated by the examples of Section~\ref{sec:applications}.

ODOS has many similarities with optimal experimental design. In both fields,  the design problem is solved by maximizing  a utility function under a parametric statistical model. The same utility functions can be used in both fields and the computational methods for optimization share many common features. From the theoretical point of view, the main difference remains that an observational design is a probability measure in the space of observational processes while an experimental design is a probability measure in the space of actual variables. It follows that in observational design each unit can be selected only once (at a given time point).  In experimental design, the design space can be continuous and the same value of a design variable can be repeated as many times as needed. It is also possible to consider hybrid studies that contain elements from both experimental and observational studies. The design space is then a product of the space of experimental variables and the space of observational processes.

Decisions on measurement times and variables to be measured need to be made in both experimental and observational studies. Some decisions in experimental design can actually be interpreted as decisions on observational design. For instance, in intervention studies with longitudinal follow-up measurements, the times of the response measurements must be decided and this decision is usually independent of  the treatment allocation.

ODOS can be criticized for its sensitivity to assumptions about specific parametric models. This criticism is not unique to ODOS but applies  to experimental design and sample size calculations as well. In the Bayesian framework, model uncertainty can be taken into account using hierarchical structures, hyperpriors and Bayesian model averaging. The statistical model used for the design optimization does not bind the hands of the analyst if a better model is found in the analysis phase. Naturally,  the design will be sub-optimal for the new model but it is still likely to be better than a non-optimized design if the analysis model and the design model share common features.

The transition from simple designs to optimized designs will emphasize the importance of careful documentation in planning and communication.   Researchers should, for instance, provide sufficient details on utility functions, prior distributions and computational methods used in ODOS. Complicated designs could be also illustrated using graphical models \citep{dagdesign}. However, this is also the challenge of ODOS.
%Complex study designs and missing data mechanisms can be illustrated using causal models with design
Often the researcher have difficulties in formulating the utility function because the planned study has multiple goals. For instance, there may be several regression parameters to be estimated and it is disputable whether D-optimality, A-optimality or some other compound criterion provides the best way to summarize the overall precision of the estimates. Approaches developed in the field of multiple-criteria decision making \citep{miettinen1999nonlinear} are potentially applicable for finding and visualizing nondominated solutions \citep{miettinen2014survey} in ODOS. In general, more research should be put also on development of utility functions.

The framework offers many interesting questions for further research. The development of efficient computational methods is one of the most important challenges. This includes both inventing  new methods and modifying optimization methods developed for optimal experimental design.  Unintended missing data provides an additional challenge for the optimization.  Finding new applications of ODOS, for instance in epidemiology, environmental sciences and marketing research, is another topic for future research. Extensive simulations with realistic parameters are needed for increasing understanding about the potential benefits of ODOS in various applications. It would be also interesting to investigate optimal designs for problems where the researcher can both observe and intervene the system.

Many of the applications presented in Section~\ref{sec:applications} were based on re-analyses of existing data. These analyses suggest that significant cost savings could be achieved if ODOS was applied in the planning of new studies. The real application of ODOS requires an open-minded principal investigator who wishes to abandon the tradition in order to try a new idea. The bottleneck in many research projects might be the need for a skilled statistician both in the planning and analysis of the study.

Compared with traditional observational designs, ODOS leads to improved precision and cost savings. On the other hand, ODOS is more complicated to implement. It is expected that the benefits will exceed the disadvantages at least in large studies and studies with expensive measurements.

%-value of information (VoI): Difference between the maximum attainable expected utility with and without executing a study that provides new observations.
%-value of control (VoC): Difference between the maximum expected utiliy with and without the ability to manipulate the state of nature.

%Could we use VoC and VoI dynamically to decide when to manipulate and when to observe: unifying framework fot experimental and observational design?

\section*{Acknowledgements}
The authors thank Antti Penttinen, Elja Arjas and Mikko Sillanp{\"a}{\"a} for useful comments. Kaisa Miettinen and Markus Hartikainen are acknowledged for discussions on multiple-criteria decision making. JK has been supported by Academy of Finland (grant numbers 266251 and 311877). JV has been supported by Academy of Finland (grant numbers 266349 and 304531) and the Research Funds of the University of Helsinki (decision No. 465/51/2014).

%The work belongs to the profiling area ``Decision analytics utilizing causal models and multiobjective optimization'' (DEMO) supported by Academy of Finland (grant number 311877).

%We will use bibtex for references
\bibliographystyle{apalike}
\bibliography{design}

\end{document}